\documentclass{article}

\usepackage[final]{neurips_2025}

\usepackage[utf8]{inputenc} 
\usepackage[T1]{fontenc}    
\usepackage{hyperref}       
\usepackage{url}            
\usepackage{booktabs}       
\usepackage{amsfonts}       
\usepackage{nicefrac}       
\usepackage{microtype}      
\usepackage{xcolor}         

\usepackage{bm}             
\usepackage{graphicx}
\usepackage{subcaption}
\usepackage{amsmath}

\usepackage{pgfplots}
\pgfplotsset{compat=1.18}
\usepackage{pgfplotstable}

\usepackage{amsmath,amsfonts,amssymb}

\usepackage{algorithm}
\usepackage{geometry}
\usepackage{graphicx}
\usepackage{booktabs}

\usepackage{algpseudocode}

\usepackage{listings}

\usepackage{tikz}
\usetikzlibrary{shapes,arrows,positioning}

\title{Convergence of finite element right-hand-side computation from finite difference data}

\author{%
  Stefan Schoder\\
  Institute of Fundamentals and Theory of Electrical Engineering\\
  Graz University of Technology\\
  8010 Graz, Austria \\
  \texttt{stefan.schoder@tugraz.at} 
}

\begin{document}

\maketitle

\begin{abstract}
This work presents two integration methods for field transfer in computational aeroacoustics and in coupled field problems, using the finite element method to solve the acoustic field. Firstly, a high-order Gaussian quadrature computes the finite element right-hand side. In contrast, the (flow) field provided by the finite difference mesh is mapped by higher-order B-Splines or a Lagrangian function. Secondly, the cut-cell or supermesh integration with geometric clipping. For each method, the accuracy, performance characteristics, and computational complexity are analyzed. As a reference, the trapezoidal integration rule was computed from the finite difference results. The high-order quadrature converges as the B-Spline interpolation order increases, and the finite difference results and mesh resolutions are consistent. The supermesh approach eliminates interpolation and approximation errors at the grid-to-mesh level and improves accuracy. This behaviour is universal for smooth or strongly oscillating field quantities, which will be shown in a comparative study between the Lighthill-like source term and the source term of the perturbed convective wave equation for subsonic flows. 
\end{abstract}

\section{Introduction}

Hybrid finite difference-finite element methods for computational aeroacoustics \cite{Schoder2019} require accurate field transfer \cite{Schoder2021} between structured finite difference grids and unstructured finite element method meshes. The finite element right-hand side assembly for acoustic source terms requires the evaluation of
\begin{equation}
\boldsymbol{b}_i = \int_{\Omega} N_i(\boldsymbol{x}) f(\boldsymbol{x}) \, d\Omega \, ,
\end{equation}
where $N_i(\boldsymbol{x})$ are finite element shape functions and $f(\boldsymbol{x})$ is the source field defined on a structured finite difference grid. Note that this approach is not limited to finite differences but can also be applied when a structured or unstructured finite volume scheme provides the input field. In the following, the implementation details, accuracy, and performance of the presented schemes are presented. A comparative study evaluated the specific differences of different aeroacoustic source terms of acoustic equations, like the ones derived in \cite{Lighthill1951,schoder_acoustic_2023,schoder_aeroacoustic_2024,schoder_perturbed_2026,Phillips1960,Lilley1974}. The visualizations in \cite{schoder_acoustic_2023,schoder_aeroacoustic_2024,schoder_perturbed_2026} show that the Lighthill-like source terms are more spatially varying and have sharper transitions as the source term exhibiting by the aeroacoustic wave equation based on Pierce operator or the extension of the perturbed convective wave equation to subsonic flows. The visual inspections suggest that interpolating smoother source terms may be less prone to errors arising from spatial discretisation.

\section{Methods}

\subsection{High-Order Gaussian Quadrature with Higher-Order Accurate Field Interpolation}

The finite element method right-hand side assembly uses Gaussian quadrature within each element
\begin{equation}
\boldsymbol{b}_i = \sum_{e=1}^{N_e} \int_{\Omega_e} N_i(\boldsymbol{x}) f_{\text{interp}}(\boldsymbol{x}) \, d\Omega \approx \sum_{e=1}^{N_e} \sum_{q=1}^{Q} w_q N_i(\boldsymbol{\xi}_q) f_{\text{interp}}(\boldsymbol{x}_q) |\boldsymbol{J}(\boldsymbol{\xi}_q)|\, ,
\end{equation}
where $Q = n_{\text{gauss}}^2$ is the number of quadrature points per element, $\boldsymbol{J}$ is the Jacobian matrix, and $f_{\text{interp}}$ is interpolated from the finite difference grid points to the Gauss quadrature points. It is important that for this type of right-hand side integration, both the interpolation order of the field function as well as the quadrature order must be sufficiently high to converge properly. Furthermore, mesh information is exchanged between the finite difference grid and the finite element mesh. The field interpolation is carried out by B-Splines supporting the orders $p \in \{1,2,3,4,5\}$ (as provided by \texttt{scipy.interpolate.RectBivariateSpline} in Python 3.10.13 scipy 1.11.4) or Lagrangian polynomials of the orders $p \in \{1,3\}$. 

For the B-spline interpolation in 2D, let \( f_{ij} = f(x_i,y_j) \) be sampled on a uniform Cartesian grid with spacing \(h\). The tensor-product B-spline interpolant of degree \(p\) is defined as
\begin{equation}
f_{\mathrm{interp}}(x,y)
=
\sum_{i}\sum_{j}
B_i^{(p)}(x)\, B_j^{(p)}(y)\, c_{ij},
\label{eq:Bspline}
\end{equation}
where \( B_i^{(p)} \) denotes the cardinal B-spline basis function of degree \(p\),
and \( c_{ij} \) are the spline coefficients obtained from the data \( f_{ij} \). The spline coefficients \( c_{ij} \) are determined by 
\begin{equation}
f(x_i,y_j)
=
\sum_{m}\sum_{n}
B_m^{(p)}(x_i)\, B_n^{(p)}(y_j)\, c_{mn},
\qquad \forall\, i,j.
\end{equation}
This yields a linear system of equations for the coefficients. 
For sufficiently smooth functions \(f \in C^{p+1}\), B-spline interpolation of
degree \(p\) achieves the error estimate
\begin{equation}
\| f - f_{\mathrm{interp}} \| = \mathcal{O}\!\left(h^{p+1}\right).
\end{equation}

The Lagrange interpolant of degree $p$ is defined as
\begin{equation}
f_{\mathrm{interp}}(x, y) =
\sum_{i}\sum_{j}
L_i^{(p)}(x) L_j^{(p)}(y) f_{ij},
\label{eq:Lagrange}
\end{equation}
where $L_i^{(p)}(x)$ denotes the Lagrange basis polynomial of degree $p$, defined as
\begin{equation}
L_i^{(p)}(x) = \prod_{\substack{k=0 \ k \neq i}}^{p} \frac{x - x_k}{x_i - x_k}.
\end{equation}
The interpolation is performed in two steps: first, a 1D Lagrange interpolation is applied along the $x$-direction to compute intermediate values; second, the values are then interpolated along the $y$-direction to obtain the final interpolated value. For $p=1$, this reduces to a bilinear interpolation. The bilinear interpolant uses two stencil points (respectively, the finite difference points) in each direction
\begin{equation}
f_{\mathrm{interp}}(x, y) =
(1 - d_\mathrm{x})(1 - d_\mathrm{y}) f_{00} +
d_\mathrm{x}(1 - d_\mathrm{y}) f_{10} +
(1 - d_\mathrm{x})d_\mathrm{y} f_{01} +
d_\mathrm{x} d_\mathrm{y} f_{11},
\label{eq:bilinear}
\end{equation}
where $d_\mathrm{x} = \frac{x - x_0}{x_1 - x_0}$ and $d_\mathrm{y} = \frac{y - y_0}{y_1 - y_0}$. The Lagrange interpolation achieves the following error estimate for sufficiently smooth functions $f \in C^{p+1}$
\begin{equation}
| f - f_{\mathrm{interp}} | = \mathcal{O}\left(h^{p+1}\right)\, .
\end{equation}

\subsubsection{Algorithm}

\begin{algorithm}[H]
\caption{Right-hand side Assembly by High-Order Quadrature}
\begin{algorithmic}[1]
\Require finite difference data $(x_{\text{fd}}, y_{\text{fd}}, f_{\text{fd}})$, finite element method mesh, interpolation order $p$, Gauss order $n_g$
\Ensure right-hand side vector $\boldsymbol{b} \in \mathbb{R}^{N_{\text{nodes}}}$
\State \textbf{Phase 1: Interpolator Construction}
\State Initialize interpolator with order $p$ and finite difference grid points
\State \textbf{Phase 2: Quadrature Point Generation}
\For{each finite element $e = 1, \ldots, N_e$}
    \State Generate $n_g \times n_g$ Gauss points in the local element reference coordinates $\boldsymbol{\xi}_q$
    \State Map to physical coordinates: $\boldsymbol{x}_q = \boldsymbol{X}(\boldsymbol{\xi}_q)$
\EndFor
\State \textbf{Phase 3: Field Interpolation}
\State Vectorized interpolation: $\boldsymbol{f}_{\text{gauss}} = \text{interpolate}(f_{\text{fd}}, \{\boldsymbol{x}_q\})$
\State \textbf{Phase 4: Assembly}
\For{each finite element $e$}
    \For{each quadrature point $q$}
        \State Evaluate shape functions: $\boldsymbol{N}_q = \boldsymbol{N}(\boldsymbol{\xi}_q)$
        \State Compute Jacobian: $J_q = \det(\boldsymbol{J}(\boldsymbol{\xi}_q))$
        \State Accumulate the right-hand side vector: $\boldsymbol{b}_{\text{local}} \leftarrow \boldsymbol{b}_{\text{local}} + w_q \boldsymbol{N}_q f_{\text{gauss}}[q] J_q$
    \EndFor
    \State Add (assemble) to the global right-hand side vector: $\boldsymbol{b} \leftarrow \boldsymbol{b}_{\text{local}}$
\EndFor
\end{algorithmic}
\end{algorithm}

\subsubsection{Convergence Analysis}

\begin{figure}
    \centering
    \includegraphics[width=0.7\linewidth,trim= 0 0 0 1.2cm,clip]{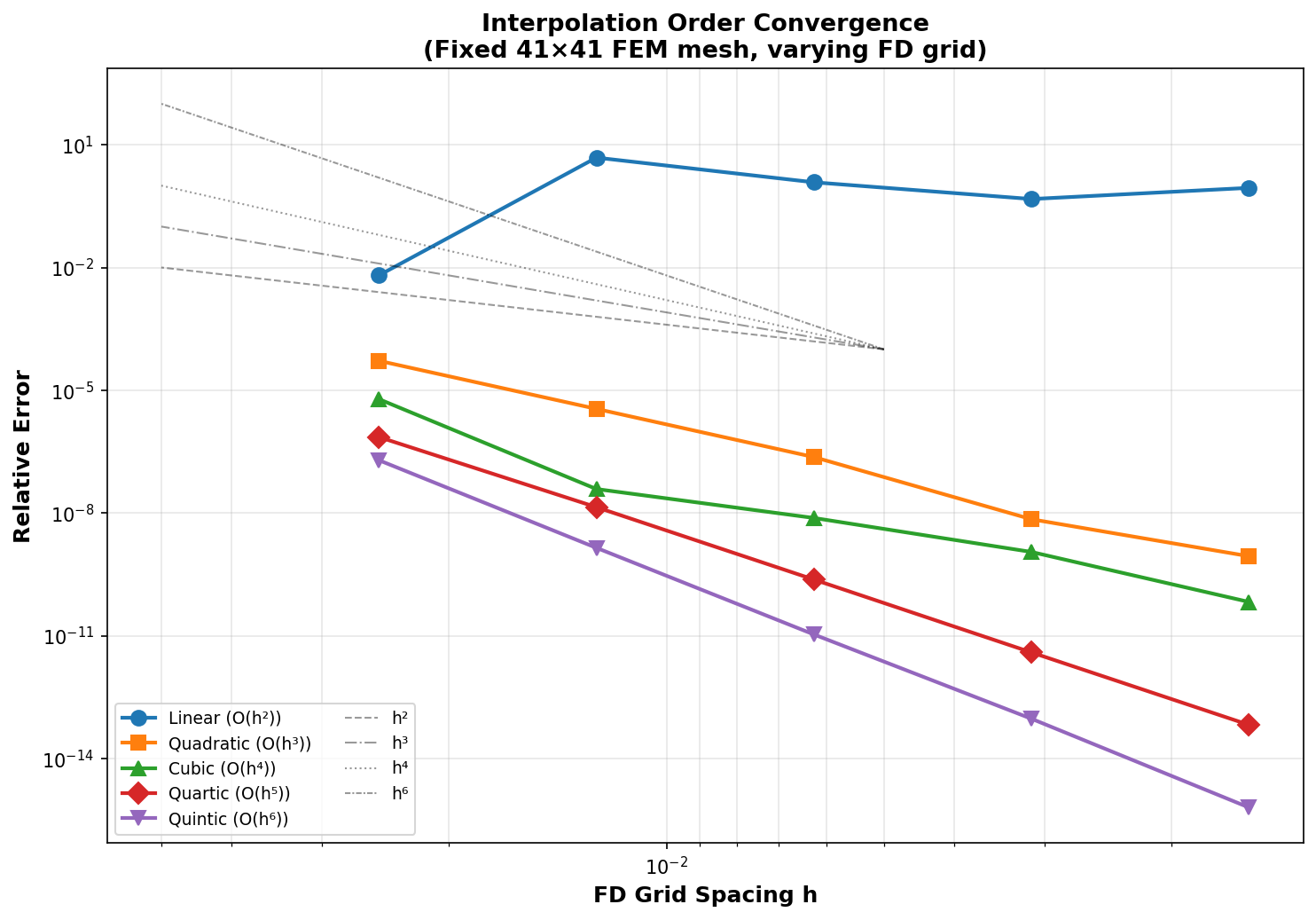}
    \caption{Interpolation error analysis as a function of finite difference grid spacing and interpolation order.}
    \label{fig:P1}
\end{figure}

\begin{figure}
    \centering
    \includegraphics[width=0.9\linewidth,trim= 0 0 0 1.6cm,clip]{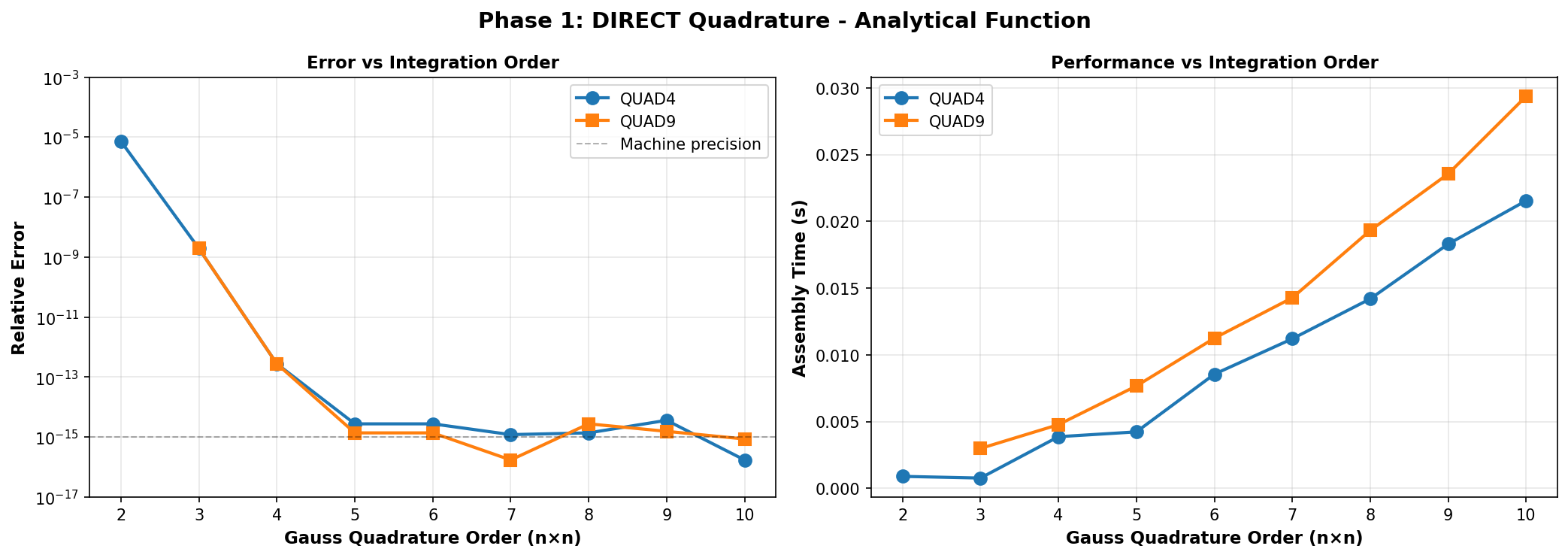}
    \caption{Interpolation error (left) and right-hand side assembly time of algorithm 1 (right) as a function of the Gauss quadrature order.}
    \label{fig:P2}
\end{figure}

\begin{figure}
    \centering
    \includegraphics[width=0.9\linewidth,trim= 0 0 0 1.6cm,clip]{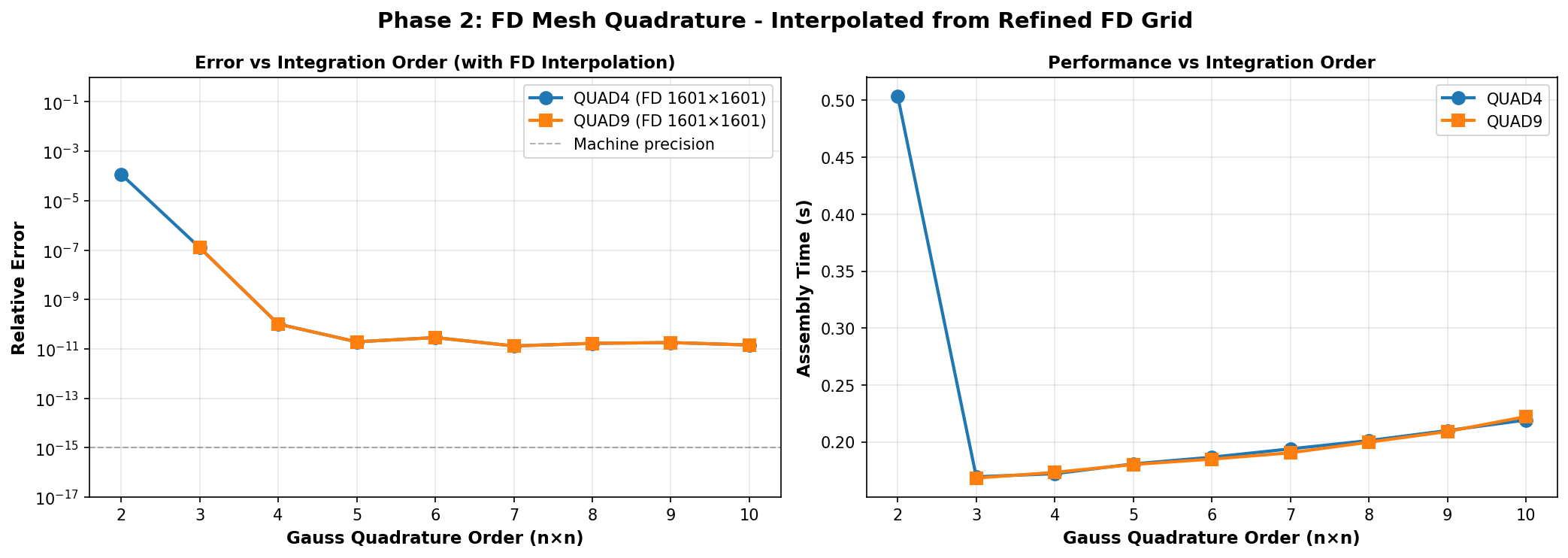}
    \caption{Interpolation error (left) and right-hand side assembly time of algorithm 1 (right) as a function of the Gauss quadrature order.}
    \label{fig:P3}
\end{figure}

The convergence analysis assesses the accuracy of the finite-difference (fd) to finite-element method (fem) data transfer pipeline, which consists of three error components
\begin{equation}
\varepsilon_{\text{total}} = \varepsilon_{\text{interp}} + \varepsilon_{\text{quad}} + \varepsilon_{\text{fem}} \, .
\end{equation}
For polynomial interpolation of order $p$, finite element basis function polynomial order $r$, and Gaussian quadrature of order $2q-1$
\begin{align}
\varepsilon_{\text{interp}} &= \mathcal{O}(h_{\text{fd}}^{p+1}) \\
\varepsilon_{\text{quad}} &= \mathcal{O}(h_{\text{fem}}^{2q}) \quad \text{if } 2q-1 \geq p \\
\varepsilon_{\text{fem}} &= \mathcal{O}(h_{\text{fem}}^{r+1}) \, .
\end{align}
$\varepsilon_{\mathrm{interp}}$ is the interpolation error from finite difference grid to finite element Gauss point locations, $\varepsilon_{\mathrm{quad}}$ is the quadrature error in the integration, and $\varepsilon_{\mathrm{fem}}$ is the finite element method spatial discretization error. The interpolation order accuracy was confirmed by studying the analytical function $f(x,y) = \sin(k_x x) \sin(k_y y)$ with $k_x = k_y = 2.5\pi$. In doing so, the grid spacing $h_{\text{fd}} \in \{2.5\cdot10^{-2}, \ldots, 1.5\cdot10^{-3}\}$ m was modified on a unit cube domain $[0;1]^2$, and the evaluation (finite element 40x40 elements) mesh was fixed and 9-point Gauss quadrature per element used. The following relative error was computed $e_{\mathrm{interp}}(h) = \left\| f_{\mathrm{interp}} - f_{\mathrm{analytic}} \right\|_{L^2(\Omega_{\mathrm{fem}})}/\left\| f_{\mathrm{analytic}} \right\|_{L^2(\Omega_{\mathrm{fem}})}$. The B-spline interpolation was used to compute the evaluation in Fig.~\ref{fig:P1}. In addition, a second study evaluated the convergence of the integral for the analytic function with $k_x = k_y = 4.5\pi$. The results presented in Fig.~\ref{fig:P2} show that the quadrature converges for a smooth function towards machine precision as the Gauss quadrature order is increased. For this test case, the analytic function was used to obtain the function values to be integrated directly, without interpolation. In a third testcase, the analytic function was first evaluated at a finite-difference location and then interpolated to the Gauss quadrature points using a Lagrangian bicubic $p=3$ interpolation. The results are provided in Fig.~\ref{fig:P3}. In the first execution, the calculation was rather slow because it was executed the first time and the respective function caches were created. In the error analysis, we can identify that for a Gauss quadrature order larger then four the interpolation error takes over and limits the overall precision of the algorithm. This is consistent with the following error estimate $E_{\mathrm{quad}}(h) = \max(h_{\text{fd}}^{p+1}, h_{\text{fem}}^{2q})$, being in the limiting case $E_{\mathrm{quad}}(h) = \max((1/1600)^{4}, (1/40)^{8}) = (1/40)^{8}$.

The practical implication of this study is that one cannot usually apply finite-difference grid refinement after the solution is already obtained. Still, if one wants to use this algorithm, high-order interpolation schemes are essential to consider to reduce the interpolation error that arises when increasing the quadrature order of the right-hand-side integration. A rule of thumb says both must be increased consistently to reduce the overall error. As a reference for the integral, the finite difference trapezoidal rule presented in the appendix can be used to quantify the error on the finite element mesh relative to the integral computed on the finite difference grid. This procedure is strongly recommended, since otherwise the algorithm cannot guarantee sufficient convergence in practice, where the field quantities may not be "sufficiently" smooth. Therefore, we are seeking to integrate aeroacoustic sources into a method that conservatively accounts for the right-hand-side source term (i.e., the complete energy from the flow field is transferred to the acoustic mesh without reduction) globally and locally \cite{Schoder2021}. The method presented in \cite{Schoder2021} only considers constant field values within one finite volume cell. For partial differential equation satisfaction while super resolving or interpolating the results, the hybrid Least Squares Support Vector Regression can be used \cite{babaei_super-resolution_2025}.

\subsection{Cut-Cell-based Supermesh Integration}
In computational physics simulations, transferring fields across incompatible discretisations is challenging, especially when the integral of the field quantity (energy exchange) must be conserved. The cut-cell or supermesh (geometric overlay of two non-conforming meshes) integration method provides an integral-preserving field transfer by computing geometric intersections and performing quadrature on the intersected meshes, resulting in cut polygons. The problem statement is that, we have given a finite difference (or finite volume) grid $\Omega_{\text{fd}} = \{(x_i, y_j) : i=0,\ldots,n_x-1, \; j=0,\ldots,n_y-1\}$, with given finite difference field values $f_{\text{fd}}(x, y)$ defined on finite difference grid, and we have given a finite element mesh $\Omega_{\text{fem}}$ with $N_e$ elements and $N_n$ nodes. The goal is to compute the finite element right-hand side (load) vector 
\begin{equation}
    \boldsymbol{b}_i = \int_{\Omega_{\text{fem}}} N_i(\boldsymbol{x}) \, f_{\text{fd}}(\boldsymbol{x}) \, d\Omega \, .
\end{equation}
The challenges in computing this integral are that the discretisation is often non-conforming, finite difference cells and finite elements have arbitrary overlap, the geometric complexity of computing polygon intersections is high, and the computational cost is therefore substantial compared to ordinary field interpolation methods.

\subsubsection{Domain Decomposition}
Let $\Omega_{\text{fem}}^e$ denote the $e$-th finite element and $\Omega_{\text{fd}}^{(i,j)}$ the finite difference cell at indices $(i,j)$, we can reformulate the integral to 
\begin{equation}
    \int_{\Omega_{\text{fem}}^e} N_k(\boldsymbol{x}) \, f_{\text{fd}}(\boldsymbol{x}) \, d\Omega = 
    \sum_{(i,j)} \int_{\Omega_{\text{fem}}^e \cap \Omega_{\text{fd}}^{(i,j)}} N_k(\boldsymbol{x}) \, f_{\text{fd,interp}}(\boldsymbol{x}) \, d\Omega
\end{equation}
where $\mathcal{P}_{e,(i,j)} = \Omega_{\text{fem}}^e \cap \Omega_{\text{fd}}^{(i,j)}$ defines an intersection polygon. The $f_{\text{fd,interp}}$ can be computed using the B-Spline \eqref{eq:Bspline}, Lagrangian \eqref{eq:Lagrange}, or bilinear interpolation \eqref{eq:bilinear}.

\subsubsection{Integration}
The integral is evaluated using Gauss quadrature at Gauss points distributed over the physical subcell of the supermesh. This does not mean that a parametric space exists for each subcell. The $f_{\text{fd,interp}}(\boldsymbol{x})$ can be evaluated on the real physical coordinates. The integration occurs in the reference coordinate system of the target element, not the subcell. Therefore, the mapping from the real physical coordinates (Gauss points) to the reference coordinates of the target element must be inverted. 

The integration domain $\mathcal{P}_{e,(i,j)} \subset \Omega_{\text{fem}}^e $ is determined by geometric intersections in physical space, rather than fixed reference coordinates.
Consequently, the quadrature points $\boldsymbol{x}_q$ must be placed to capture the shape of the physical cut region $\mathcal{P}_{e,(i,j)}$.
Evaluating the basis functions $N_i$ at these physical points requires the inverse mapping $\boldsymbol{\xi}_q = \boldsymbol{X}^{-1}(\boldsymbol{x}_q)$.
For higher-order polynomial isoparametric elements or higher-order geometrically mapped finite elements, this inverse is not analytic and Newton-Raphson might obtain inversion for every quadrature point:
\begin{equation}
\int_{\mathcal{P}_{e,(i,j)}} f_{\text{fd,interp}}(\boldsymbol{x})\,N_i(\boldsymbol{x})\,\mathrm{d}\boldsymbol{x}
=
\sum_{q \in \text{cut}}
f_{\text{fd,interp}}(\boldsymbol{x}_q)\,
N_i\bigl(\boldsymbol{X}^{-1}(\boldsymbol{x}_q)\bigr)\,
\det\!\left(\frac{\partial \boldsymbol{X}}{\partial \boldsymbol{\xi}}\Bigl(\boldsymbol{X}^{-1}(\boldsymbol{x}_q)\Bigr)\right)\,
w_q.
\end{equation}
This introduces high computational cost and complexity compared to standard assembly, as the mapping $\boldsymbol{X}(\boldsymbol{\xi}_q)$ is designed to be evaluated explicitly only in the forward direction. To compute the quadrature on generalized intersection polygon $\mathcal{P}_{e,(i,j)}$, the polygon is tessellated into triangles $\{\mathcal{T}_1, \mathcal{T}_2, \ldots, \mathcal{T}_{n_T}\}$. On each triangle $\mathcal{T}_m$, the 6-point Gaussian quadrature is used to integrate the right-hand side partly. 
The finite element method right-hand side vector component on element $e$ for the local degree of freedom $k$ is calculated as 
\begin{equation}
    b_k = \sum_{(i,j)} \sum_{m=1}^{n_T^{e,(i,j)}} |\mathcal{T}_m| \sum_{q=1}^{6} w_q \, N_k(\boldsymbol{x}_q) \, f_{\text{fd,interp}}(\boldsymbol{x}_q) \, ,
\end{equation}
where $(i,j)$ are the finite difference cells intersecting element $e$, $n_T^{e,(i,j)}$ is the number of triangles in the intersection polygon $\mathcal{P}_{e,(i,j)}$, and $|\mathcal{T}_m|$ is the area of triangle $m$.

In standard finite element assembly, the right-hand side corresponding to an analytic forcing function
$f(\boldsymbol{x})$ is assembled as
\begin{equation}
\int_{\Omega_e} f(\boldsymbol{x})\,N_i(\boldsymbol{x})\,\mathrm{d}\boldsymbol{x}
=
\sum_q
f\bigl(\boldsymbol{X}(\boldsymbol{\xi}_q)\bigr)\,
N_i(\boldsymbol{\xi}_q)\,
\det\!\left(\frac{\partial \boldsymbol{X}}{\partial \boldsymbol{\xi}}(\boldsymbol{\xi}_q)\right)\,
w_q,
\end{equation}
without requiring any inverse mapping or nonlinear iteration.

\subsubsection{Physical space to target element reference mapping}
A smooth, invertible transformation
\begin{equation}
\boldsymbol{x} = \boldsymbol{X}(\boldsymbol{\xi}),
\qquad
\boldsymbol{\xi} \in \hat{\Omega},\ \boldsymbol{x} \in \Omega,
\end{equation}
where $\hat{\Omega}$ denotes the reference element and $\boldsymbol{X}$ is typically a polynomial or spline-based mapping of degree $p \ge 1$. This can be an expansion using the same basis functions as the finite element basis functions for isoparametric finite elements, or a different mapping function.
Given a physical point $\boldsymbol{x}^\star \in \Omega$ (the seeded Gauss points in real physical space), a pointwise evaluation of a
finite element geometry mapping requires the corresponding reference coordinates $\boldsymbol{\xi}^\star$ satisfying
\begin{equation}
\boldsymbol{X}(\boldsymbol{\xi}^\star) = \boldsymbol{x}^\star.
\end{equation}
For linear mappings this inversion is analytic. However, for curvilinear or high-order geometry mapping $\boldsymbol{X}$ is nonlinear and does not admit a closed-form inverse. Consequently, the reference coordinates $\boldsymbol{\xi}^\star$ must be obtained by solving the nonlinear system
\begin{equation}
\boldsymbol{F}(\boldsymbol{\xi}) = \boldsymbol{X}(\boldsymbol{\xi}) - \boldsymbol{x}^\star = \boldsymbol{0}\, .
\end{equation}
This can be performed using a Newton-Raphson iteration, which is provided in the appendix.

\subsubsection{Geometry Processing}
Axis-Aligned Bounding Box (AABB) indexing is used to find finite difference cells that intersect a finite element efficiently. We use binary search on the structured grid to get a list of potential finite difference cell candidate indices $(i,j)$ inside a target element bounding box $[x_{\min}, x_{\max}] \times [y_{\min}, y_{\max}]$. In general, the clipping of a finite element against an AABB is performed by the Sutherland-Hodgman polygon clipping algorithm \cite{sutherland_reentrant_1974} or the CGAL library.

\subsubsection{Algorithm}
The algorithm is separated into a setup phase and an execution phase, which is executed at every timestep of the simulation data.
\begin{algorithm}[H]
\caption{Supermesh Integration Pipeline}
\begin{algorithmic}[1]
\Require finite difference grid $(x_{\text{fd}}, y_{\text{fd}}, f_{\text{fd}})$, finite element method mesh, tolerance $\varepsilon = 10^{-10}$
\Ensure right-hand side vector $\boldsymbol{b}$

\State \textbf{Phase 1: AABB Spatial Indexing}
\For{each finite element $e$}
    \State Compute element bounding box using $\text{AABB}_e = [\min x_e, \min y_e, \max x_e, \max y_e]$
    \State Find candidate finite difference cells $\mathcal{C}_e$
\EndFor

\State \textbf{Phase 2: Geometric Clipping}
\For{each element $e$ and cell $c \in \mathcal{C}_e$}
    \State Create polygons $P_e = \text{Polygon}(e)$, $P_c = \text{Polygon}(c)$
    \State $\Omega_{e,c} = \text{Intersection}(P_e, P_c)$
\EndFor

\State \textbf{Phase 3: Sub-Tessellation}
\For{each intersection polygon $\Omega_{e,c}$}
    \State Triangulate $\{\mathcal{T}_k\} = \text{Tesselate}(\Omega_{e,c})$
    \For{each triangle $\mathcal{T}_k$}
        \State Map six-point Gauss quadrature rule to $\mathcal{T}_k$: $\{(\boldsymbol{x}_{kq}, w_{kq})\}_{q=1}^6$
    \EndFor
\EndFor

\State \textbf{Phase 4: Field Evaluation}
\For{each Gauss point $\boldsymbol{x}_{kq}$}
    \State \textbf{finite difference Evaluation:} $f_{\text{fd,interp}}(\boldsymbol{x}_{kq}) = \text{Interp}(\boldsymbol{x}_{kq}, c)$
    \State \textbf{finite element method Evaluation:}
    \State \quad Newton-Raphson: $\boldsymbol{\xi} = \text{InverseMap}(\boldsymbol{x}_{kq}, e, \varepsilon)$
    \State \quad Shape functions: $\boldsymbol{N}(\boldsymbol{\xi})$, Jacobian: $J(\boldsymbol{\xi})$
    \State \textbf{Integration and Assembly:} $\boldsymbol{b} \leftarrow \boldsymbol{b} +  w_{kq} \boldsymbol{N}(\boldsymbol{\xi}) f_{\text{fd}}(\boldsymbol{x}_{kq}) J(\boldsymbol{\xi})$
\EndFor
\end{algorithmic}
\end{algorithm}

\section{Comparative study}
\begin{figure}[htbp]
    \centering
    \begin{subfigure}[b]{0.48\textwidth}
        \includegraphics[width=\textwidth]{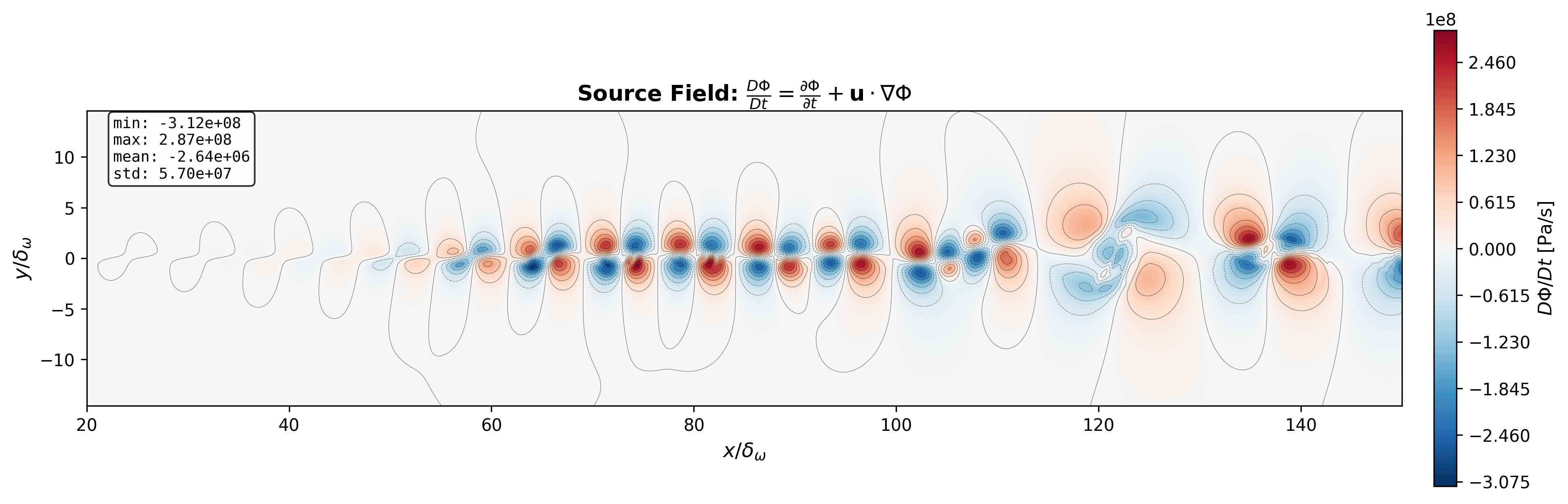}
        \caption{$D\Phi/Dt$}
    \end{subfigure}
    \hfill
    \begin{subfigure}[b]{0.48\textwidth}
        \includegraphics[width=\textwidth]{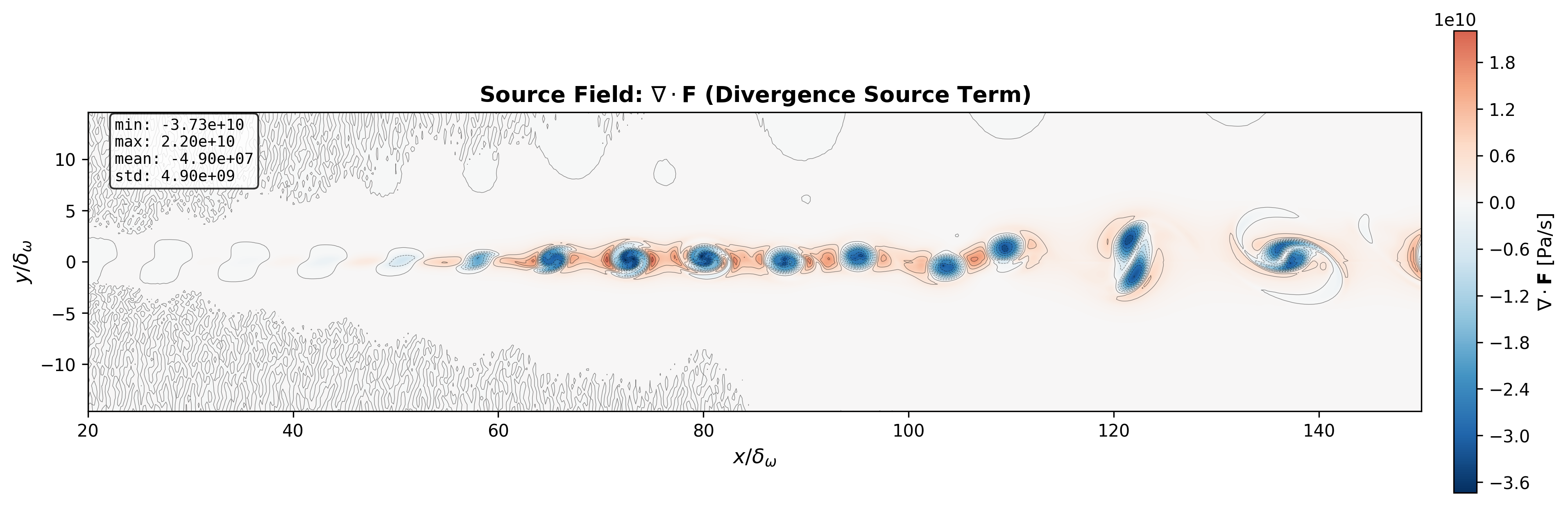}
        \caption{$\nabla \cdot \mathbf{F}$}
    \end{subfigure}
    \caption{Source fields used for integration.}
    \label{fig:source_fields}
\end{figure}

\begin{figure}[htbp]
    \centering
    \includegraphics[width=0.9\textwidth,trim= 0 0 0 1.4cm,clip]{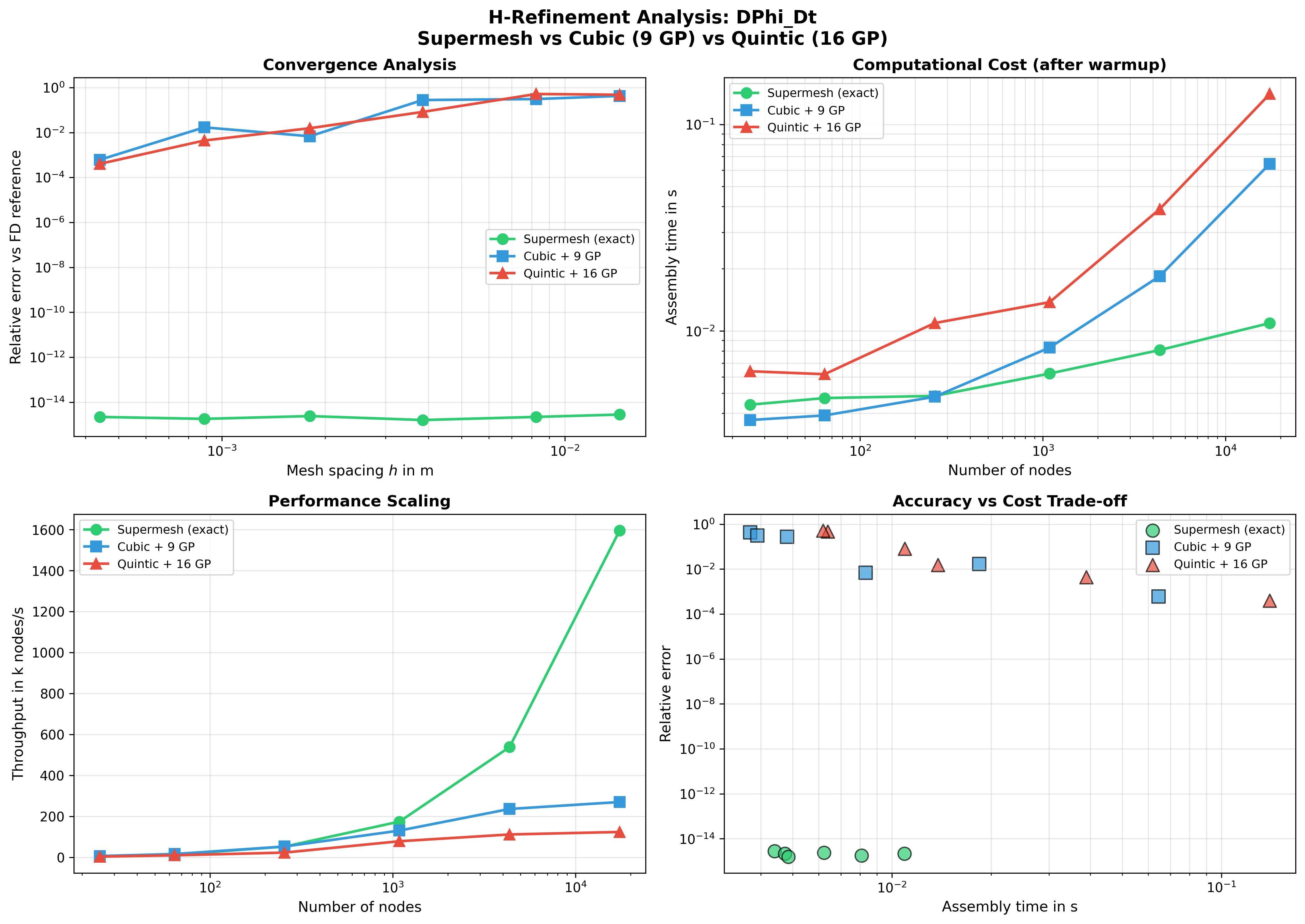}
    \caption{H-refinement analysis with CFD data for the material derivative $D\Phi_p/Dt$. Comparison between supermesh integration, cubic interpolation, and quintic interpolation.}
    \label{fig:hrefinement_DPhi}
\end{figure}

\begin{figure}[htbp]
    \centering
    \includegraphics[width=0.9\textwidth,trim= 0 0 0 1.4cm,clip]{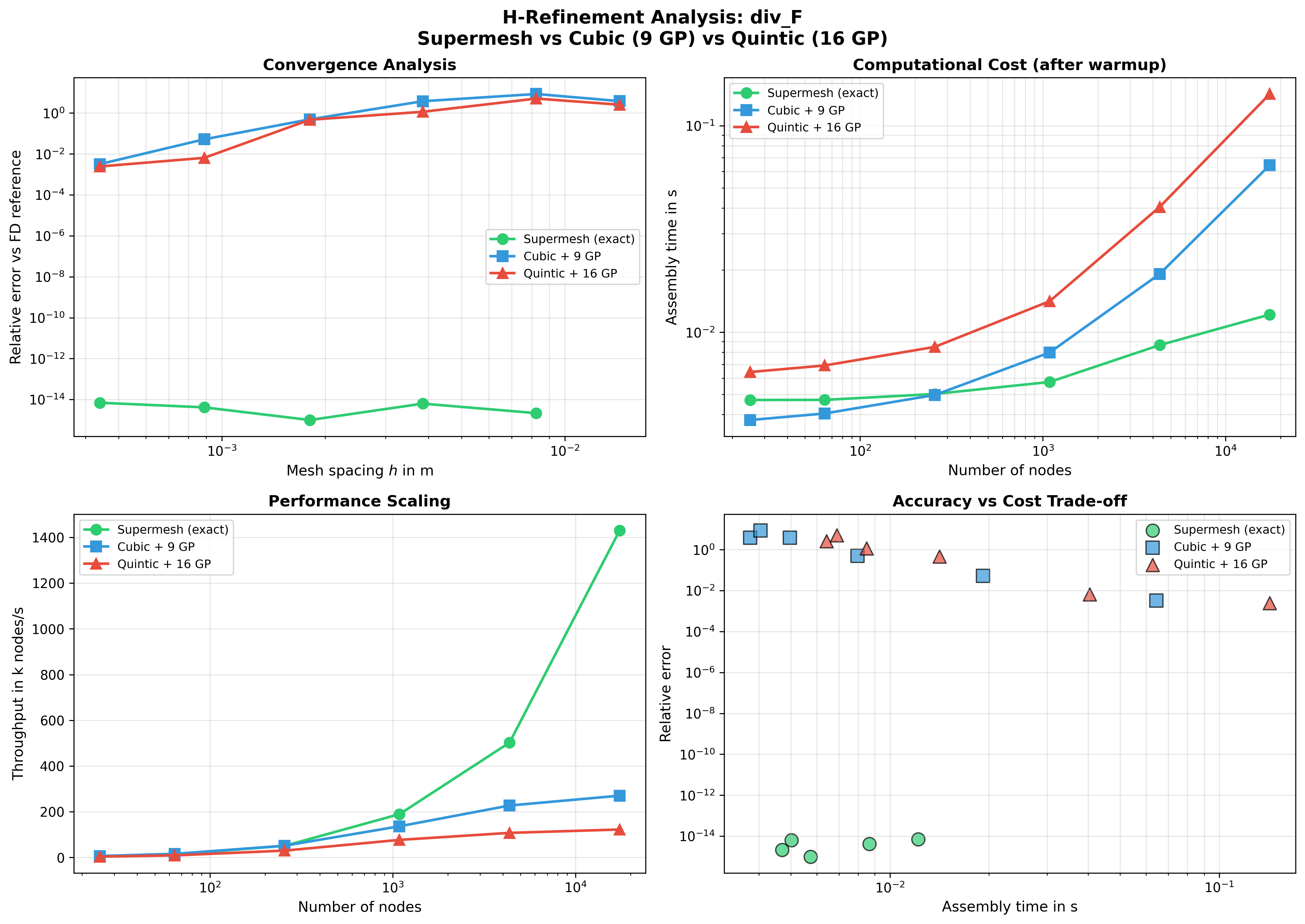}
    \caption{H-refinement analysis with CFD data for the source divergence $\nabla \cdot \mathbf{F}$. Comparison between supermesh integration, cubic interpolation, and quintic interpolation.}
    \label{fig:hrefinement_divF}
\end{figure}

The following CFD data from the case \cite{schoder_perturbed_2026,vincent_application_2023,schoder_aeroacoustic_2022} is transformed from the finite difference grid $x/\delta_\omega \in [20, 150]$, $y/\delta_\omega \in [-15, 15]$ to the finite element mesh. $D\Phi_p/Dt$ and $\nabla \cdot \mathbf{F} = -\nabla \cdot \left[ 
    (\mathbf{u}_v \cdot \nabla)\mathbf{u}_v + 
    (\mathbf{u}_0 \cdot \nabla)\mathbf{u}_v + (\mathbf{u}_v \cdot \nabla)\mathbf{u}_0 
    \right]$ as an exemplary source term are considered (see Fig~\ref{fig:source_fields}), for more detail we refer to \cite{schoder_perturbed_2026}. In figure \ref{fig:hrefinement_DPhi} and \ref{fig:hrefinement_divF}, the cut-cell supermesh geometric clipping, the bicubic B-spline interpolation with 3×3 Gauss quadrature, and the quintic B-spline interpolation with 4×4 Gauss quadrature are shown. All timings were measured after a warm-up run to ensure JIT compilation and cache-warming effects were excluded. The supermesh method converges to machine-precision for all mesh resolutions. In contrast, interpolation-based methods exhibit systematic errors that converge with mesh refinement, improved quadrature, and higher interpolation order. Furthermore, it can be seen that the method is already highly optimized and outperforms ordinary interpolation methods for many elements to be processed; however, the standalone interpolation methods with quadrature are not optimized at this point. There is a systematic difference between the source terms; the $D\Phi_p/Dt$ shows less error for the quadrature-based interpolation methods.
This analysis (see Fig.~\ref{fig:quadrature_analysis}) investigates the effect of quadrature order of the interpolation-based method on the accuracy of source term integration using interpolation-based B-Spline methods. The mesh resolution is fixed at $\Delta x/\delta_\omega = \Delta y/\delta_\omega = 1$. Both cubic and quintic interpolation converge with increasing quadrature order but do not reach machine precision. The interpolation error for higher-order quadrature rules is dominated by the B-spline reconstruction of the finite difference field, not the quadrature error. 
In Figure \ref{fig:weak_scaling_realdata}, the weak scaling of the source term integration is shown, with no notable difference between the source terms.

\section{Conclusions}

\begin{table}[htbp]
\centering
\caption{Comparison of integration methods for source term transformation from a finite difference grid to a finite element mesh.}
\label{tab:performance_comparison}
\begin{tabular}{@{}llll@{}}
\toprule
 & \textbf{Supermesh} & \textbf{B-Spline}  \\ \midrule
Conservation Error & Machine precision & Interpolation/Quadrature limited  \\
Target $h$-Refinement & Exact at all scales & Systematic error floor (interpolation) \\
Primary error source & None (Geometric) & Reconstruction order of the interpolation \\
Mapping required & Inverse ($X^{-1}$) & Forward ($X$) \\
Scaling & Linear (Weak) & Linear (Weak) \\
Source term robustness & Provided & Smoother Sources Preferred \\
\bottomrule
\end{tabular}
\end{table}

\section*{Acknowledgments}
The authors acknowledge the support of rAI. Parts of the methods are implemented in openCFS \cite{Schoder2022openCFS}, openCFS-data\cite{schoder_opencfs-data_2023}, and will be available in pyCFS-data \cite{wurzinger2024pycfs} soon.

{
\small

\bibliographystyle{unsrt}
\bibliography{zotero}

@misc{schoder_opencfs-data_2023,
	title = {{openCFS}-{Data}: {Data} {Pre}-{Post}-{Processing} {Tool} for {openCFS}},
	copyright = {arXiv.org perpetual, non-exclusive license},
	shorttitle = {{openCFS}-{Data}},
	url = {https://arxiv.org/abs/2302.03637},
	doi = {10.48550/ARXIV.2302.03637},
	urldate = {2026-01-19},
	publisher = {arXiv},
	author = {Schoder, Stefan and Roppert, Klaus},
	year = {2023},
	note = {Version Number: 3},
}

@misc{babaei_super-resolution_2025,
	title = {Super-{Resolution} of {Elliptic} {PDE} {Solutions} {Using} {Least} {Squares} {Support} {Vector} {Regression}},
	url = {http://arxiv.org/abs/2512.09967},
	doi = {10.48550/arXiv.2512.09967},
	urldate = {2026-01-19},
	publisher = {arXiv},
	author = {Babaei, Maryam and Rucz, Peter and Kaltenbacher, Manfred and Schoder, Stefan},
	month = dec,
	year = {2025},
	note = {arXiv:2512.09967 [math]},
}

@article{sutherland_reentrant_1974,
	title = {Reentrant polygon clipping},
	volume = {17},
	issn = {0001-0782, 1557-7317},
	url = {https://dl.acm.org/doi/10.1145/360767.360802},
	doi = {10.1145/360767.360802},
	language = {en},
	number = {1},
	urldate = {2026-01-18},
	journal = {Communications of the ACM},
	author = {Sutherland, Ivan E. and Hodgman, Gary W.},
	month = jan,
	year = {1974},
	pages = {32--42},
}

@article{schoder_perturbed_2026,
	title = {Perturbed convective wave equation for low-to-medium {Mach} number subsonic flows},
	volume = {623},
	issn = {0022460X},
	url = {https://linkinghub.elsevier.com/retrieve/pii/S0022460X25006224},
	doi = {10.1016/j.jsv.2025.119549},
	language = {en},
	urldate = {2026-01-17},
	journal = {Journal of Sound and Vibration},
	author = {Schoder, S. and Bagheri, E. and Bogey, C. and Bailly, C.},
	month = feb,
	year = {2026},
	pages = {119549},
}

@misc{Schoder2022openCFS,
	title = {{openCFS}: {Open} source finite element software for coupled field simulation – part acoustics},
	doi = {10.48550/ARXIV.2207.04443},
	publisher = {arXiv},
	author = {Schoder, S. and Roppert, K.},
	year = {2022},
	note = {tex.archiveprefix: arxiv
tex.copyright: Creative Commons Attribution 4.0 International},
}

@article{wurzinger2024pycfs,
	title = {{pyCFS}-data: {Data} processing framework in python for {openCFS}},
	journal = {arXiv preprint arXiv:2405.03437},
	author = {Wurzinger, Andreas and Schoder, Stefan},
	year = {2024},
}

@article{vincent_application_2023,
	title = {Application of the complex differentiation method to the sensitivity analysis of aerodynamic noise},
	volume = {264},
	issn = {00457930},
	url = {https://linkinghub.elsevier.com/retrieve/pii/S0045793023001901},
	doi = {10.1016/j.compfluid.2023.105965},
	language = {en},
	urldate = {2025-10-14},
	journal = {Computers \& Fluids},
	author = {Vincent, Hugo and Bogey, Christophe},
	month = oct,
	year = {2023},
	pages = {105965},
}

@inproceedings{schoder_aeroacoustic_2022,
	address = {Southampton, UK},
	title = {Aeroacoustic wave equation based on {Pierce}'s operator applied to the sound generated by a mixing layer},
	isbn = {978-1-62410-664-4},
	url = {https://arc.aiaa.org/doi/10.2514/6.2022-2896},
	doi = {10.2514/6.2022-2896},
	language = {en},
	urldate = {2025-10-03},
	booktitle = {28th {AIAA}/{CEAS} {Aeroacoustics} 2022 {Conference}},
	publisher = {American Institute of Aeronautics and Astronautics},
	author = {Schoder, Stefan and Kaltenbacher, Manfred and Spieser, Étienne and Vincent, Hugo and Bogey, Christophe and Bailly, Christophe},
	month = jun,
	year = {2022},
}

@article{schoder_aeroacoustic_2024,
	title = {Aeroacoustic {Source} {Potential} {Based} on {Poisson}’s {Equation}},
	volume = {62},
	issn = {0001-1452, 1533-385X},
	url = {https://arc.aiaa.org/doi/10.2514/1.J063792},
	doi = {10.2514/1.J063792},
	language = {en},
	number = {7},
	urldate = {2025-10-03},
	journal = {AIAA Journal},
	author = {Schoder, Stefan and Bagheri, Eman and Spieser, Étienne},
	month = jul,
	year = {2024},
	pages = {2772--2782},
}

@article{schoder_acoustic_2023,
	title = {Acoustic {Modeling} {Using} the {Aeroacoustic} {Wave} {Equation} {Based} on {Pierce}’s {Operator}},
	volume = {61},
	issn = {0001-1452, 1533-385X},
	url = {https://arc.aiaa.org/doi/10.2514/1.J062558},
	doi = {10.2514/1.J062558},
	language = {en},
	number = {9},
	urldate = {2025-09-04},
	journal = {AIAA Journal},
	author = {Schoder, Stefan and Spieser, Étienne and Vincent, Hugo and Bogey, Christophe and Bailly, Christophe},
	month = sep,
	year = {2023},
	pages = {4008--4017},
}

@article{Schoder2021,
	title = {Application limits of conservative source interpolation methods using a low {Mach} number hybrid aeroacoustic workflow},
	volume = {29},
	doi = {10.1142/S2591728520500322},
	number = {1},
	journal = {Journal of Theoretical and Computational Acoustics},
	publisher = {World Scientific},
	author = {Schoder, S. and Wurzinger, A. and Junger, C. and Weitz, M. and Freidhager, C. and Roppert, K. and Kaltenbacher, M.},
	year = {2021},
	pages = {2050032},
}

@article{Lighthill1951,
	title = {On sound generated aerodynamically {I}. {General} theory},
	volume = {211},
	doi = {10.1098/rspa.1952.0060},
	journal = {Proceedings of the Royal Society of London},
	author = {{M. J. Lighthill}},
	year = {1952},
	pages = {564--587},
}

@techreport{Lilley1974,
	title = {On the noise from jets},
	institution = {AGARD CP-131},
	author = {{G. M. Lilley}},
	year = {1974},
}

@article{Phillips1960,
	title = {On the generation of sound by supersonic turbulent shear layers},
	volume = {9},
	doi = {10.1017/S0022112060000888},
	journal = {Journal of Fluid Mechanics},
	author = {{O. M. Phillips}},
	year = {1960},
	pages = {1--18},
}

@article{Schoder2019,
	title = {Hybrid aeroacoustic computations: {State} of art and new achievements},
	volume = {27},
	doi = {10.1142/S2591728519500208},
	number = {04},
	journal = {Journal of Theoretical and Computational Acoustics},
	publisher = {World Scientific},
	author = {Schoder, S. and Kaltenbacher, M.},
	year = {2019},
	pages = {1950020},
}
}

\appendix

\section{Reference integration by the trapezoidal rule}

The 2D trapezoidal rule provides reference integration on structured finite difference grids
\begin{equation}
I_{\text{trap}} = \int_{\Omega} f(x,y) \, dx \, dy \approx \sum_{i=0}^{n_x-1} \sum_{j=0}^{n_y-1} w_{ij} f_{ij}
\end{equation}
where weights $w_{ij}$ account for grid point weighting and non-uniform spacing. For smooth functions $f \in C^2(\Omega)$, the trapezoidal rule accuracy is
\begin{equation}
\varepsilon_{\text{trap}} = I_{\text{exact}} - I_{\text{trap}} = \mathcal{O}(h^2)
\end{equation}
where $h = \max(\Delta x_{\max}, \Delta y_{\max})$ is the maximum grid spacing. The weights are computed as follows for a structured grid. For grid coordinates $\{x_i\}_{i=0}^{n_x-1}$ and $\{y_j\}_{j=0}^{n_y-1}$ with spacings $\Delta x_i = x_{i+1} - x_i$ and $\Delta y_j = y_{j+1} - y_j$:

\begin{align}
w_{ij} = \begin{cases}
\frac{1}{4} \Delta x_0 \Delta y_0 & \text{corner: } (0,0) \\
\frac{1}{4} \Delta x_{n_x-2} \Delta y_0 & \text{corner: } (n_x-1,0) \\
\frac{1}{4} \Delta x_0 \Delta y_{n_y-2} & \text{corner: } (0,n_y-1) \\
\frac{1}{4} \Delta x_{n_x-2} \Delta y_{n_y-2} & \text{corner: } (n_x-1,n_y-1) \\
\frac{1}{2} \Delta x_0 \cdot \frac{1}{2}(\Delta y_{j-1} + \Delta y_j) & \text{left edge: } i=0, j \in [1,n_y-2] \\
\frac{1}{2} \Delta x_{n_x-2} \cdot \frac{1}{2}(\Delta y_{j-1} + \Delta y_j) & \text{right edge: } i=n_x-1, j \in [1,n_y-2] \\
\frac{1}{2}(\Delta x_{i-1} + \Delta x_i) \cdot \frac{1}{2} \Delta y_0 & \text{bottom edge: } i \in [1,n_x-2], j=0 \\
\frac{1}{2}(\Delta x_{i-1} + \Delta x_i) \cdot \frac{1}{2} \Delta y_{n_y-2} & \text{top edge: } i \in [1,n_x-2], j=n_y-1 \\
\frac{1}{2}(\Delta x_{i-1} + \Delta x_i) \cdot \frac{1}{2}(\Delta y_{j-1} + \Delta y_j) & \text{interior: } i,j \in [1,n-2]
\end{cases}
\end{align}

\subsection*{Algorithm Workflow}

\begin{algorithm}[H]
\caption{2D Trapezoidal Integration}
\begin{algorithmic}[1]
\Require $x_{\text{fd}} \in \mathbb{R}^{n_x}$, $y_{\text{fd}} \in \mathbb{R}^{n_y}$, $f_{\text{fd}} \in \mathbb{R}^{n_y \times n_x}$
\Ensure $I_{\text{trap}} \in \mathbb{R}$
\State Compute spacings: $\Delta x_i = x_{i+1} - x_i$, $\Delta y_j = y_{j+1} - y_j$
\State Initialize weight matrix: $W \in \mathbb{R}^{n_y \times n_x}$
\For{$j = 1$ to $n_y-2$}
    \For{$i = 1$ to $n_x-2$}
        \State $W[j,i] = \frac{1}{2}(\Delta x_{i-1} + \Delta x_i) \cdot \frac{1}{2}(\Delta y_{j-1} + \Delta y_j)$
    \EndFor
\EndFor
\State Compute edge weights
\State Compute corner weights
\Return $I_{\text{trap}} = \sum_{i,j} W[j,i] \cdot f_{\text{fd}}[j,i]$
\end{algorithmic}
\end{algorithm}

\section{Target element vectorized Newton-Raphson}
Using isoparametric bilinear quadrilateral finite elements in 2d, the residual of the inverse mapping function is given by
\begin{equation}
\boldsymbol{R}(\boldsymbol{\xi}) = \boldsymbol{X}(\boldsymbol{\xi}) - \boldsymbol{x}^* = \sum_{j=1}^4 N_j(\boldsymbol{\xi}) \boldsymbol{x}_j - \boldsymbol{x}^* = \boldsymbol{0}
\end{equation}
where, the solution $\boldsymbol{\xi}^*$ (when the residual converges to e.g. $\|\boldsymbol{R}\| < 10^{-10}$) can be obtained Newton-Raphson iteration
\begin{equation}
\boldsymbol{\xi}^{(k+1)} = \boldsymbol{\xi}^{(k)} - \boldsymbol{J}^{-1}(\boldsymbol{\xi}^{(k)}) \boldsymbol{R}(\boldsymbol{\xi}^{(k)}) \, ,
\end{equation}
where the gradient of the residual is defined by the Jacobian $\boldsymbol{J}(\boldsymbol{\xi}) = \frac{\partial \boldsymbol{R}}{\partial \boldsymbol{\xi}} = \frac{\partial \boldsymbol{X}}{\partial \boldsymbol{\xi}}$.

This can be efficiently implemented using a vectorized Newton-Raphson scheme that processes all quadrature points for an element $j$-index summation. Use the physical coordinates $\boldsymbol{x}_{in}= [x_1, \ldots, x_n]^T$, $\boldsymbol{y}_{in} = [y_1, \ldots, y_n]^T$ with a batch number of $n$ points. 

\begin{align}
    \boldsymbol{F}_x^{(k)} &= \sum_{j=1}^{4} N_j(\boldsymbol{\xi}_r^{(k)}, \boldsymbol{\eta}_r^{(k)}) x_j - \boldsymbol{x}_{in} \\
    \boldsymbol{F}_y^{(k)} &= \sum_{j=1}^{4} N_j(\boldsymbol{\xi}_r^{(k)}, \boldsymbol{\eta}_r^{(k)}) y_j - \boldsymbol{y}_{in} \\
    \boldsymbol{J}_{11}^{(k)} &= \sum_{j=1}^{4} \frac{\partial N_j}{\partial \xi_r}(\boldsymbol{\xi}_r^{(k)}, \boldsymbol{\eta}_r^{(k)}) x_j \\
    \boldsymbol{J}_{12}^{(k)} &= \sum_{j=1}^{4} \frac{\partial N_j}{\partial \eta_r}(\boldsymbol{\xi}_r^{(k)}, \boldsymbol{\eta}_r^{(k)}) x_j \\
    \boldsymbol{J}_{21}^{(k)} &= \sum_{j=1}^{4} \frac{\partial N_j}{\partial \xi_r}(\boldsymbol{\xi}_r^{(k)}, \boldsymbol{\eta}_r^{(k)}) y_j \\
    \boldsymbol{J}_{22}^{(k)} &= \sum_{j=1}^{4} \frac{\partial N_j}{\partial \eta_r}(\boldsymbol{\xi}_r^{(k)}, \boldsymbol{\eta}_r^{(k)}) y_j \\
    \det(\boldsymbol{J})^{(k)} &= \boldsymbol{J}_{11}^{(k)} \boldsymbol{J}_{22}^{(k)} - \boldsymbol{J}_{12}^{(k)} \boldsymbol{J}_{21}^{(k)} \\
    \boldsymbol{\xi}_r^{(k+1)} &= \boldsymbol{\xi}_r^{(k)} - \frac{\boldsymbol{J}_{22}^{(k)} \boldsymbol{F}_x^{(k)} - \boldsymbol{J}_{12}^{(k)} \boldsymbol{F}_y^{(k)}}{\det(\boldsymbol{J})^{(k)}} \\
    \boldsymbol{\eta}_r^{(k+1)} &= \boldsymbol{\eta}_r^{(k)} - \frac{\boldsymbol{J}_{11}^{(k)} \boldsymbol{F}_y^{(k)} - \boldsymbol{J}_{21}^{(k)} \boldsymbol{F}_x^{(k)}}{\det(\boldsymbol{J})^{(k)}}
\end{align}

\newpage
\section{Weak Scaling}

\begin{figure}[htbp]
    \centering
    \begin{subfigure}[b]{0.48\textwidth}
        \includegraphics[width=\textwidth,trim= 0 0 0 1.4cm,clip]{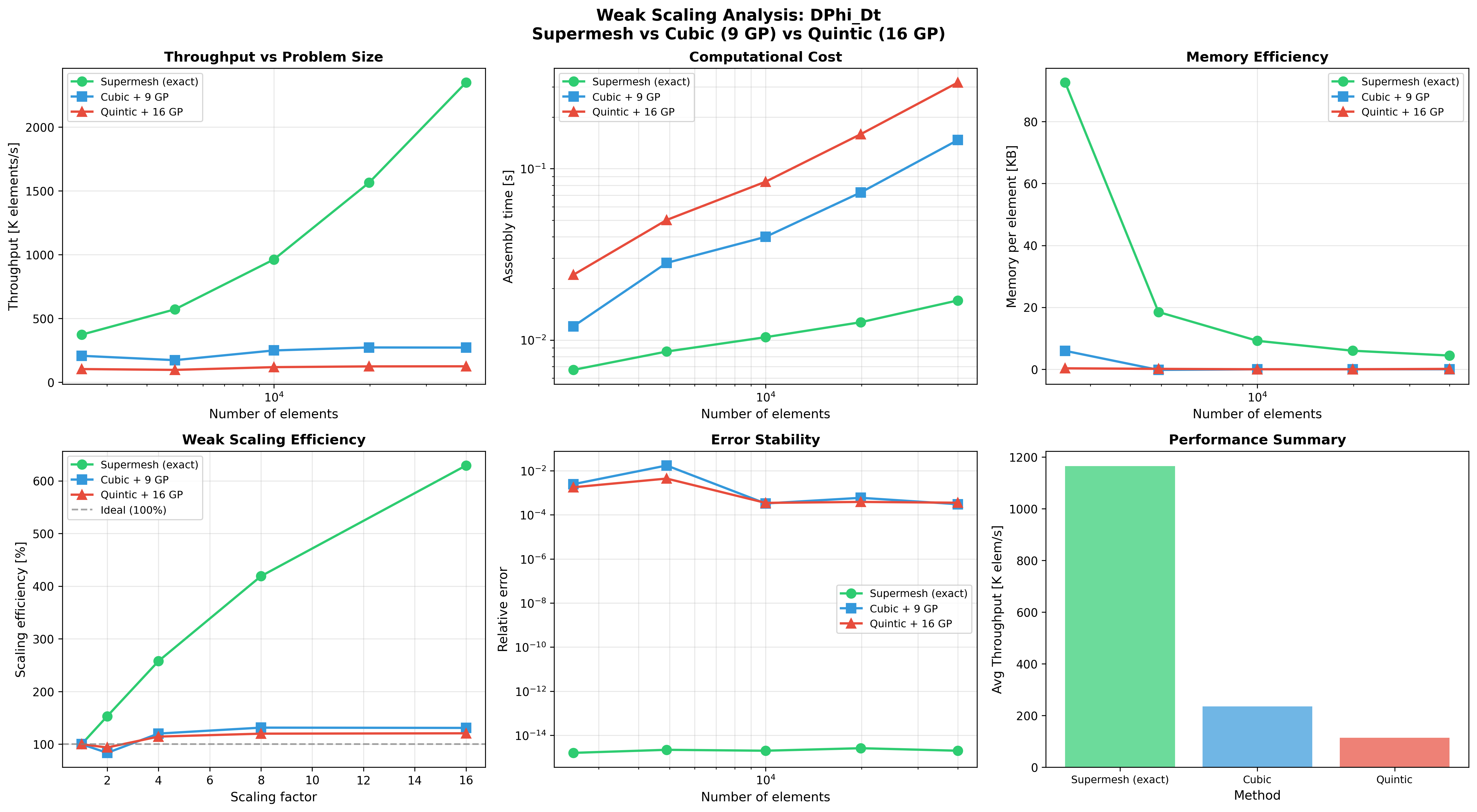}
        \caption{$D\Phi/Dt$}
    \end{subfigure}
    \hfill
    \begin{subfigure}[b]{0.48\textwidth}
        \includegraphics[width=\textwidth,trim= 0 0 0 1.4cm,clip]{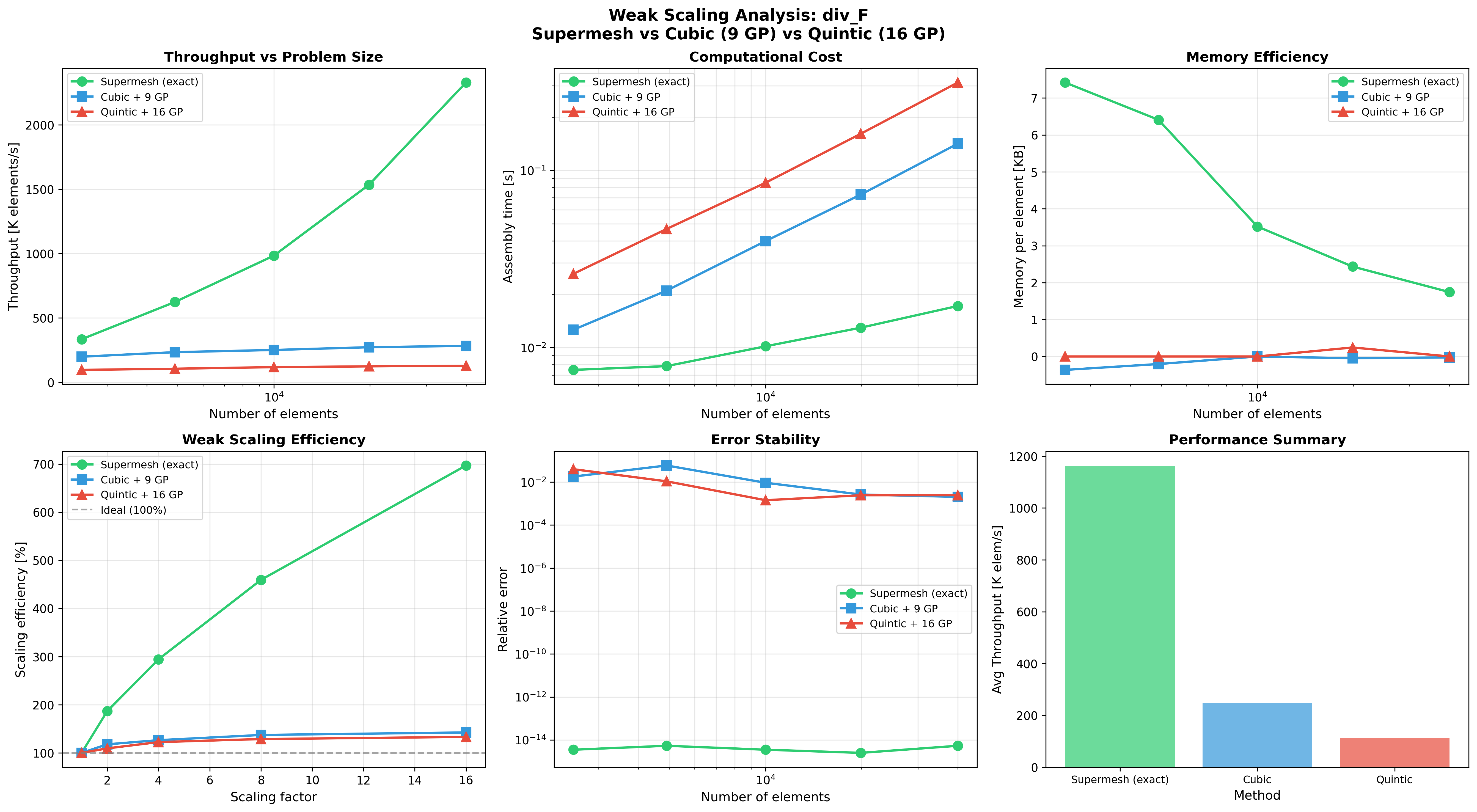}
        \caption{$\nabla \cdot \mathbf{F}$}
    \end{subfigure}
    \caption{Weak scaling analysis with CFD data.}
    \label{fig:weak_scaling_realdata}
\end{figure}

\newpage
\section{Interpolation Study}

\begin{figure}[htbp]
    \centering
    \begin{subfigure}[b]{0.98\textwidth}
        \includegraphics[width=\textwidth,trim= 0 0 0 1.4cm,clip]{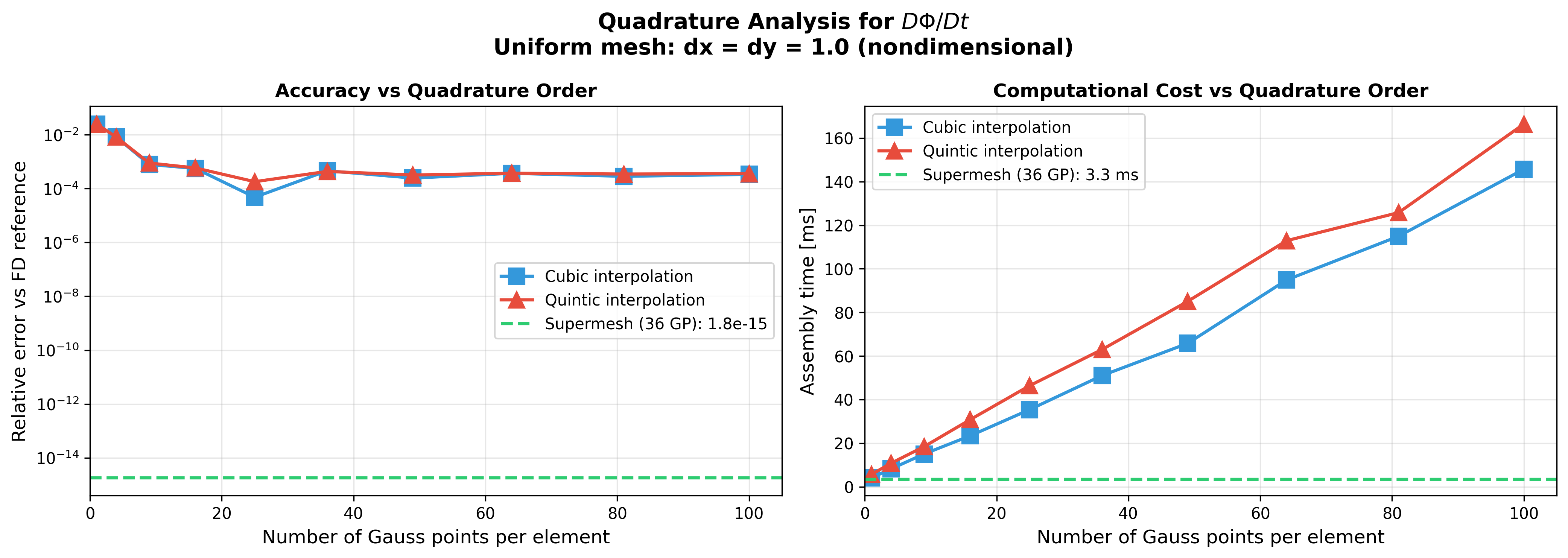}
        \caption{$D\Phi/Dt$}
    \end{subfigure}
    \hfill
    \begin{subfigure}[b]{0.98\textwidth}
        \includegraphics[width=\textwidth,trim= 0 0 0 1.4cm,clip]{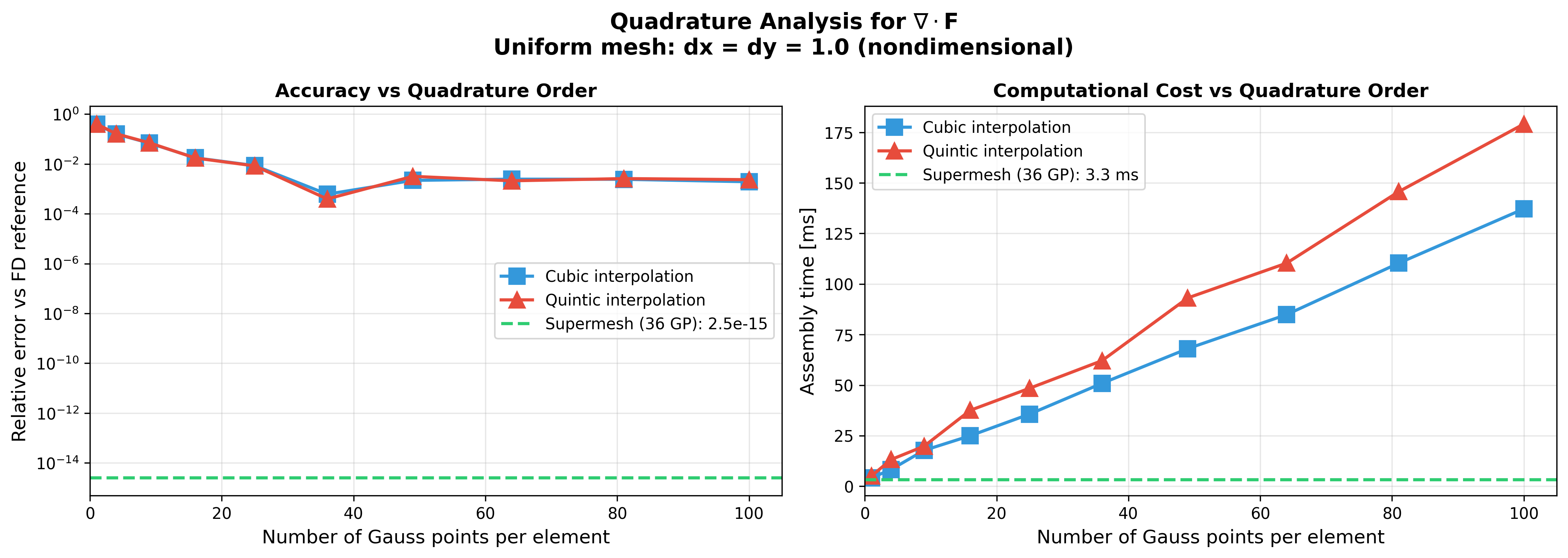}
        \caption{$\nabla \cdot \mathbf{F}$}
    \end{subfigure}
    \caption{Effect of quadrature order on integration accuracy. The horizontal dashed line indicates the supermesh reference accuracy.}
    \label{fig:quadrature_analysis}
\end{figure}

\end{document}